\documentclass{amsart}
\usepackage{stmaryrd}
\usepackage[english]{babel}
\usepackage{amsfonts,amssymb}
\newtheorem{theorem}{Theorem}[section]
\newtheorem*{theorem*}{Theorem}
\newtheorem{corollary}[theorem]{Corollary}
\newtheorem{lemma}[theorem]{Lemma}

\theoremstyle{definition}
\newtheorem{definition}[theorem]{Definition}
\theoremstyle{remark}
\newtheorem{remark}[theorem]{Remark}
\newtheorem{remarks}[theorem]{Remarks}

\def\cprime{$'$}
\newcommand{\N}{\mathbf N}
\newcommand{\Z}{\mathbf Z}
\newcommand{\A}{\mathcal A}
\newcommand{\B}{\mathcal B}
\def\epsi{\varepsilon}
\def\ti{-\allowhyphens}
\def\nimply{\Longarrownot\Longrightarrow}
\def\se{\subseteq}
\def\bsl{\backslash}
\begin{document}
\title[Amenable actions]{Amenable actions, free products\\ and a fixed point property}
\author{Yair Glasner and Nicolas Monod}

\address{Y.G.: University of Illinois at Chicago, Department of Mathematics, Chicago, IL 60607 USA}
\address{N.M.: University of Chicago, Department of Mathematics, Chicago, IL 60637 USA}
\begin{abstract}
We investigate the class of groups admitting an action on a set with an invariant mean. It turns out that many free products admit such an action. We give a complete characterisation of such free products in terms of a strong fixed point property.
\end{abstract}
\maketitle
\section{Introduction}
\subsection{}
In the early 20th century, the construction of Lebesgue's measure was followed by the discovery of the Banach-Hausdorff-Tarski paradoxes (\cite{Lebesgue01}, \cite{Hausdorff14}, \cite{Banach23}, \cite{Tarski29}; see also~\cite{Harpe_Leb}). This prompted von Neumann~\cite{vonNeumann29} to study the following general question:

\smallskip

Given a group $G$ acting on a set $X$, when is there an \emph{invariant mean} on $X\,$? 

\begin{definition}
An invariant mean is a $G$\ti invariant map $\mu$ from the collection of subsets of $X$ to $[0,1]$ such that (i)~$\mu(A\cup B) = \mu(A) + \mu(B)$ when $A\cap B=\varnothing$ and (ii)~$\mu(X)=1$. If such a mean exists, the action is called \emph{amenable}.
\end{definition}

\begin{remarks}
(1)~For the study of the classical paradoxes, one also considers normalisations other than~(ii). (2)~The notion of amenability later introduced by Zimmer~\cite{Zimmer84} is different, being in a sense dual to the above.
\end{remarks}

\subsection{}
The thrust of von Neumann's article was to show that the paradoxes, or lack thereof, originate in the structure of the group rather than the set $X$. He therefore proposed the study of \emph{amenable groups} (then ``messbare Gruppen''), i.e. groups whose action on themselves by multiplication is amenable. This direction of research turned out to be most fruitful, with influences on combinatorial group theory, ergodic theory, rigidity and semi-simple groups, harmonic analysis, operator algebras, etc.

\smallskip

However, the original question remained largely unanswered (compare Greenleaf~\cite{GreenleafBook}, \emph{Problem} p.~18 and Pier~\cite{Pier} p.~307). Whilst it is easy to see that any action of an amenable group is amenable, the converse a priori holds only for \emph{free} actions, where it is essentially tautological. Besides free actions, von Neumann describes only one other example (pp.~82--83 in~\cite{vonNeumann29}) which still almost contains a free action of a free group, noting about the general case that ``its somewhat complicated character might be found annoying''.

\smallskip

When investigating the general question of the amenability of a $G$\ti action on $X$, a few restrictions are in order (compare Greenleaf, \emph{loc. cit.}). First, one should assume the action \emph{faithful}, since otherwise one is really investigating a quotient group of $G$. Next, it is natural to consider \emph{transitive} actions, since otherwise $X$ could contain any action (e.g. a fixed point, providing an obvious invariant mean) as long as one adds a free orbit for the sake of faithfulness. Thus, we shall focus in this note on the class of all countable\footnote{\emph{countable} shall always mean of cardinality~$\aleph_0$, thus infinite.} groups that admit a faithful transitive amenable action:
$$\A = \big\{\, G \text{ countable} \ : \ G {\text{ admits a faithful transitive amenable action }}\big\}.$$
The only obvious examples of groups in~$\A$ are amenable groups (since the $G$\ti action on itself is free and transitive). Given that in the classical paradoxes the non-amenability was caused by the presence of a non-Abelian free group, the following posthumous result of E.K.~van Douwen is at first sight surprising:

\begin{theorem*}[van Douwen~\cite{vanDouwen}]
 Finitely generated non-Abelian free groups are in~$\A$.
\end{theorem*}

\subsection{}
It is easy to verify that if $G$ has \emph{Kazhdan's property~(T)}, then any amenable $G$\ti action has a finite orbit (Lemma~\ref{lem:T}). Thus for instance $\mathbf{SL}_3(\Z)$ is not in~$\A$. We propose the following:

\begin{definition}
A countable group $G$ has the \emph{fixed point property~(F)} if any amenable $G$\ti action (on a countable set) has a fixed point.
\end{definition}

Such a group is never in~$\A$; examples include infinite simple Kazhdan groups. However, $G$ need not be Kazhdan because property~(F) is preserved under finite free products. (See Section~\ref{sec:little} for details.) We say that $G$ has \emph{virtually~(F)} if a finite index subgroup of $G$ has property~(F).

\begin{remark}
The virtual property is much stronger than requiring that $G$ have a finite orbit, since it implies for instance that $G$ has a \emph{minimal} finite index subgroup. (Indeed a group with property~(F) cannot have any non-trivial finite quotient since that would provide an amenable $G$\ti set without fixed point.)
\end{remark}

\subsection{}
The main result of this note is that $G*H$ is \emph{always} in~$\A$ unless the obvious obstruction occurs:

\begin{theorem}\label{thm:main}
Let $G,H$ be any countable groups. Then $G*H\in\A$ unless $G$ has property~(F) and $H$ has virtually~(F). (Upon possibly exchanging $G$ and $H$.)

\medskip
\noindent
Moreover, this provides a necessary and sufficient characterisation.
\end{theorem}

For example, we see that for any countable $G,H$, the free product $G*H$ is in~$\A$ as soon as one of the groups is either residually finite or non-finitely-generated or amenable. Thus for instance $G*\mathbf{SL}_3(\Z)$ is in~$\A$ for any countable group $G$. Furthermore, Theorem~\ref{thm:main} leads to the following dichotomy:

\begin{corollary}\label{cor:infinite}
Let $G=*_{i=1}^n G_i$ be any free product of $2\leq n \leq \infty$ countable groups. Then either $G\in \A$ or $G$ has virtually~(F).

\medskip

\noindent
Moreover, the latter occurs if and only if $n\neq \infty$, all $G_i$ with $i>1$ have property~(F) and $G_1$ virtually has~(F). (Upon possibly reordering the factors.)
\end{corollary}

Incidentally, this shows that a group in~$\A$ can be the directed union of groups with property~(F).

\subsection{}
We summarize below some structural properties of the class~$\A$. Most of these properties are either elementary or follow easily from known results (and from Theorem~\ref{thm:main}). Statement~(\ref{P:semidirect}) provides an interesting contrast to Theorem~\ref{thm:main}. See Section~\ref{sec:little} for definitions and details.

\begin{theorem}\label{thm:properties}
For any countable groups $G, H$, the following hold:\nobreak
\begin{enumerate}
\item\label{P:amenable}
$G$ amenable $\Longrightarrow$ $G\in\A$.
\smallskip
\item\label{P:product}
$G,H\in\A$ $\Longleftrightarrow$ $G\times H\in\A$.
\smallskip
\item\label{P:*}
$G, H\in\A$ $\Longrightarrow$ $G*H\in\A$.\\
However, $G*H\in\A$ $\nimply$ $G, H\in\A$.
\smallskip
\item\label{P:semidirect}
$G, H\in\A$ $\nimply$ $G\rtimes H\in\A$, even if $G={\bf Z}^2$.
\smallskip
\item\label{P:embed}
Any countable group embeds into a group in~$\A$.
\smallskip
\item\label{P:co-amenable}
Assume $H< G$ is co-amenable. Then $H\in\A$ $\Longrightarrow$ $G\in\A$.\\
However, $G\in\A$ $\nimply$ $H\in\A$, even if $G=H\rtimes{\bf Z}$.
\smallskip
\item\label{P:T}
$G$ has Kazhdan's property~(T) $\Longrightarrow$ $G\notin\A$.
\smallskip
\item\label{P:rel_T}
Let $H\triangleleft G$ be a normal subgroup that is not of finite exponent. If the pair $(H,G)$ has the relative property~(T), then $G\notin\A$.
\smallskip
\item\label{P:F_and_T}
R.~Thompson's group $F$ is in~$\A$; non-amenable Tarski monsters are not.
\end{enumerate}
\end{theorem}

\begin{remark}\label{rem:co-amenable}
About~(\ref{P:co-amenable}): A subgroup $H< G$ is called \emph{co-amenable} if the $G$\ti action on $X= G/H$ is amenable. Thus we could rephrase all questions about~$\A$ by studying groups $G$ that admit some co-amenable subgroup $H$ such that the intersection of all conjugates of $H$ is trivial. (See~\cite{Eymard72}, \cite{Monod-Popa} for more on co-amenability.) 
\end{remark}

\begin{remark}
The starting point of this note was our observation that one can give a very short proof of van Douwen's result that $\Z*\Z\in\A$: If $\sigma$ is a transitive permutation of a countable set $X$, then any \emph{generic} choice (in Baire's sense) of a permutation $\tau$ defines a faithful transitive \emph{amenable} action of a free group with $\sigma,\tau$ as free generators. The idea that generic transformations generate a free group has been used e.g. in~\cite{Epstein71},\cite{Dixon90},\cite{Abert-Virag},\cite{Abert_laws}. Another simple proof of van Douwen's result was communicated to us by R.I.~Grigorchuk, to appear~\cite{Grigorchuk-Nekrashevych}. Using generic permutations, one can further establish that~$\A$ is closed under free product.
\end{remark}

\subsection{}
The organisation of this paper is as follows. Section~\ref{sec:lemmata} gathers basic facts about amenable actions. Section~\ref{sec:proof} is concerned with the proof of Theorem~\ref{thm:main}. Section~\ref{sec:little} supplies the proofs of the remaining statements or examples of this Introduction.

\medskip

It is a pleasure to thank R.I.~Grigorchuk for bringing to our attention the problem formulated by Greenleaf~\cite[p.~18]{GreenleafBook}. The work of the second author was partially supported by NSF grant DMS~0204601.

\section{Generalities}
\label{sec:lemmata}
\subsection{}
A \emph{$G$\ti set} is a countable set endowed with an action of the countable group $G$; a \emph{$G$\ti map} is a $G$\ti equivariant map between $G$\ti sets. Unless otherwise stated, $G$ itself is considered as a $G$\ti set under left multiplication. The group of bijections of $X$ is denoted by $X!$. By functoriality one has:

\begin{lemma}
\label{lemma:map}%
Let $X\to Y$ be a $G$\ti map of $G$\ti sets. If the $G$\ti action on $X$ is amenable, then so is its action on $Y$.\hfill\qedsymbol
\end{lemma}

\noindent
This shows notably that any action of an amenable group is amenable. In the anti-functorial direction, one checks:

\begin{lemma}\label{lem:anti}
Let $X$ be a $G$\ti set with an invariant mean $\mu$ and let $Y\se X$ be a $G$\ti invariant subset. If $\mu(Y)\neq0$, then $\mu/\mu(Y)$ yields an invariant mean on $Y$.\hfill\qedsymbol
\end{lemma}

\subsection{}
Recall that a subgroup $H<G$ is \emph{co-amenable} if the $G$\ti action on $G/H$ is amenable (e.g. if $H$ has finite index in $G$). This is equivalent to the following relative fixed point property (see~\cite{Eymard72}):

\emph{Every continuous affine $G$\ti action on a convex compact subset of a locally convex space with an $H$\ti fixed point has a $G$\ti fixed point.}

\smallskip

Applying this to the space of means on a $G$\ti set $X$, one deduces:

\begin{lemma}\label{lemma:co-amen}
Let $X$ be a $G$\ti set and $H$ a co-amenable subgroup. If the $H$\ti action on $X$ is amenable, then so is the $G$\ti action.\hfill\qedsymbol
\end{lemma}

\subsection{}
Let $H$ be a countable group, $L<H$ a subgroup and $Z$ an $L$\ti set. If $Z$ were a coset space $Z=L/K$ for some subgroup $K<L$, one would obtain a related $H$\ti set $X$ by setting $X=H/K$. This construction can be generalized to the arbitrary $L$-set $Z$ as follows.

\begin{definition}
The \emph{induced} $H$\ti set is the quotient $X=L\bsl(Z\times H)$ of $Z\times H$ by the diagonal $L$\ti action; the $H$\ti action on itself by \emph{right} multiplication (by the inverse elements) turns $X$ into an $H$\ti set.
\end{definition}

It is straightforward to verify the following.

\begin{lemma}\label{lemma:ind0}
If the $L$\ti action on $Z$ is faithful, transitive or free then the $H$\ti action on $X$ has the corresponding property. The converse holds for the latter two properties but not for faithfulness.\hfill\qedsymbol
\end{lemma}

\begin{lemma}\label{lemma:ind}
Suppose that the $L$\ti action on $Z$ is amenable. Then the $H$\ti action on $X$ is amenable if and only if $L$ is co-amenable in $H$.
\end{lemma}

\begin{remark}
However, even if $Z$ is of the form $L/K$ for $K$ normal in $L$, it can happen that the $H$\ti action on $X$ is amenable and $L<H$ is a co-amenable whilst the $L$\ti action on $Z$ is not amenable; see~\cite{Monod-Popa}.
\end{remark}

\begin{proof}[Proof of Lemma~\ref{lemma:ind}]
The map $Z\to Z\times\{e\}\se Z\times H$ descends to an $L$\ti map $Z\to X$. Therefore, by Lemma~\ref{lemma:map}, the $L$\ti action on $X$ is amenable. Thus, if $L$ is co-amenable in $H$, the $H$\ti action on $X$ is amenable by Lemma~\ref{lemma:co-amen}. The converse follows from Lemma~\ref{lemma:map} since there is a canonical $H$\ti map $X\to H/L$.
\end{proof}

\begin{lemma}\label{lem:sub:F}
Suppose $L\lhd H$ is normal and co-amenable. Then $L$ has property~(F) if and only if every amenable $H$\ti set has an $L$\ti fixed point.
\end{lemma}

\begin{proof}
Necessity is obvious. Conversely, let $Z$ be an amenable $L$\ti set; $L$ fixes a point in the induced $H$\ti set $X$ by Lemma~\ref{lemma:ind}. If this point is represented by $(z, h)$ in $L\bsl(Z\times H)$, then $z$ is fixed by $h L h^{-1} = L$.
\end{proof}

\subsection{}
The following characterisation originating in F{\o}lner's work~\cite{Folner} is well-known; see~\cite{Rosenblatt73} for a proof.

\begin{theorem}\label{thm:Folner}
A $G$\ti action on a set $X$ is amenable if and only if for any finite subset $S\se G$ and any $\epsi>0$ there exists a finite subset $A\se X$ such that
$$|A \vartriangle sA| < \epsi |A| \quad \forall s \in S.$$
\end{theorem}

\noindent
We call such a set an \emph{$(S,\epsi)$\ti F{\o}lner set}.

\begin{remark}\label{rem:orbit}
The above inequality is additive with respect to decomposing $A$ along the partition of $X$ into $G$\ti orbits. Therefore, given $S$ and $\epsi$, we can find an $(S,\epsi)$\ti F{\o}lner set contained in a single $G$\ti orbit.
\end{remark}

Since we consider the case where $G$ is countable, it follows  from Theorem~\ref{thm:Folner} that the action is amenable if and only if there exists a sequence $\{A_n\}_{n=1}^\infty$ of finite non-empty sets $A_n\se X$ such that for every $s\in G$ one has $\lim_{n\to \infty} |A_n \vartriangle s A_n|/|A_n|=0$.

\begin{definition}
A sequence as above is called a \emph{F{\o}lner sequence} for the $G$\ti action on $X$.
\end{definition}

\begin{remark}\label{rem:gen}
It suffices to check $\lim_{n\to \infty} |A_n \vartriangle s A_n|/|A_n|=0$ for all $s$ in some set generating $G$.
\end{remark}

\begin{lemma}\label{lem:big}
Let $H$ be a countable group and assume that $H$ does not have virtually~(F).

Then there exists an $H$\ti set $Y$ and a F{\o}lner sequence $\{A_n\}_{n\in\N}$ such that (i)~each $A_n$ is contained in a single $H$\ti orbit $H y_n \se Y$ and (ii)~the cardinality $|H y_n|$ converges to infinity (in $\N\cup\{\aleph_0\}$). Moreover we may assume that each $H$\ti orbit in $Y$ contains at most one set $A_n$.
\end{lemma}

\begin{proof}
Let $H=\{h_n : n\in\N\}$ be an enumeration of $H$. We claim that for every $n\in\N$ there is an $H$\ti set $Y_n$ and a finite subset $A_n\se Y_n$ contained in a single $H$\ti orbit $H y_n\se Y_n$ such that
$$\text{(1) } \forall\, i\leq n:\ |A_n \vartriangle h_i A_n| < \frac{1}{n} |A_n|, \quad\quad\quad\text{and (2) } |H y_n|>n.$$
This implies the statement of the lemma upon considering $Y=\bigsqcup_{n\in \N} Y_n$. Thus we consider for a contradiction the smallest $n\in \N$ for which the claim fails. Then every amenable $G$\ti set $X$ has an orbit of size at most $n$. Indeed, by Theorem~\ref{thm:Folner} $X$ contains a finite subset satisfying~(1) and we may assume that it lies in a single $H$\ti orbit by Remark~\ref{rem:orbit}; therefore~(2) has to fail. Considering the special case where $X$ is a disjoint union of finite quotients of $H$, we conclude that $H$ has a \emph{minimal} finite index subgroup $L$. Let now $Z$ be any amenable $L$\ti set. Since the induced $H$\ti set $X$ is amenable by Lemma~\ref{lemma:ind}, it contains a finite $H$\ti orbit by the previous discussion. But then $Z$ has a finite $L$\ti orbit and hence a fixed point since $L$ has no finite index proper subgroup. Thus $L$ has property~(F), a contradiction.
\end{proof}

\subsection{}
An idea of Kazhdan~\cite{Kazhdan67} yields:

\begin{lemma}\label{lem:K}
Let $G$ be a countable group. If $G$ is not finitely generated, then $G$ has an amenable action without finite orbits.
\end{lemma}

\begin{proof}
Let $X$ be the disjoint union of all coset spaces $G/H$, where $H$ ranges over the family of finitely generated subgroups of $G$. Then $X$ is a countable set with a natural $G$\ti action; there are no finite orbits since $G$ is not finitely generated. For any finite subset $S\se G$, the trivial coset $\langle S\rangle$ in $G/\langle S\rangle$ is fixed by every $s\in S$; thus, the set $A=\{\langle S\rangle\}$ is an $(S,\epsi)$\ti F{\o}lner set for any $\epsi>0$ and the action is amenable.
\end{proof}

\section{Amenable Actions of Free Products}
\label{sec:proof}
\subsection{}
We first explain why the restriction in Theorem~\ref{thm:main} is an obvious obstruction; the rest of Section~\ref{sec:proof} will be devoted to prove that this is the only obstruction.

\begin{lemma}\label{lem:converse}
Let $G$ be a group with property~(F) and $H$ with virtually~(F). Then $G*H$ has virtually~(F).
\end{lemma}

\noindent
In particular, such a group $G*H$ cannot belong to the class~$\A$.

\begin{proof}[Proof of the lemma]
Let $H_0<H$ be a finite index subgroup with property~(F); since $H_0$ is a minimal finite index subgroup, it is normal in $H$. We claim that the kernel $K$ of the canonical morphism $G*H\to H/H_0$ has property~(F). By Lemma~\ref{lem:sub:F}, it suffices to find a $K$\ti fixed point in any $G*H$\ti set $X$ with an invariant mean $\mu$. Fix coset representatives $h_1, \ldots, h_n\in H$ for $H/H_0$. Let $X^G$ be the set of $G$\ti fixed points. By property~(F), the $G$\ti action on $X\setminus X^G$ is not amenable. Therefore, $\mu(X\setminus X^G)=0$; indeed, otherwise, Lemma~\ref{lem:anti} would yield a $G$\ti invariant mean on $X\setminus X^G$. It follows that $\mu(X^G)=1$ and likewise $\mu(X^{H_0})=1$. Since $\mu$ is invariant, we deduce further $\mu(h^{-1}X^G)=1$ for all $h\in H$. Therefore, the set
$$X^{H_0}\ \cap\ \bigcap_{i=1}^n\,h_i^{-1} X^G$$
has mean one and hence is non-empty. For any $x$ in this set, the $H$\ti orbit $Hx = \{h_1 x, \ldots, h_n x\}$ is $G$\ti fixed and therefore consists of $K$\ti fixed points.
\end{proof}

The above proof also shows:

\begin{lemma}\label{lem:F:products}
If $G$ and $H$ have property~(F), then so does $G*H$.\hfill\qedsymbol
\end{lemma}

\subsection{}
We now establish a result slightly weaker than Theorem~\ref{thm:main}:

\begin{theorem}\label{thm:weaker}
Let $G,H$ be any countable groups. Then $G*H\in\A$ unless both $G$ and $H$ have virtually~(F).
\end{theorem}

\begin{proof}
We may assume that $H$ does not have virtually~(F). Let $Y$ be an $H$\ti set as in Lemma~\ref{lem:big}. Let $G=D\sqcup R$ be a partition of $G$ into two infinite sets. We may index the $H$\ti orbits of $Y$ by $R$ and write $Y=\bigsqcup_{g\in R} Y_g$. We endow the set $X=H\sqcup Y$ with the natural $H$\ti action. Given any injective map $\beta:G\to X$ we obtain a $G$\ti action on $X$ as follows: We transport by $\beta$ the action of $G$ on itself to the corresponding $G$\ti action on $\beta(G)$ and $G$ acts trivially on $X\setminus \beta(G)$. Denote by $G_\beta< X!$ the resulting subgroup; since $H$ acts faithfully on $X$ we may consider it as a subgroup $H<X!$ and thus by universal property we have a canonical epimorphism $G*H\to \langle G_\beta, H\rangle$. We shall show that for a suitable choice of $\beta$ the resulting $G*H$\ti action on $X$ is faithful, transitive and amenable.

\smallskip

We shall only consider maps $\beta$ such that $\beta(D)\se H$ and $\beta(R)\se Y$. We start by determining $\beta$ on $R$. Let $\{A_n\}_{n\in \N}$ be a sequence of subsets of $Y$ as in Lemma~\ref{lem:big}. For each $g\in R$ we make any choice $\beta(g)\in Y_g$. Since $Y_g$ meets (and then contains) at most one $A_n$, we may require that $\beta(g)\notin A_n$ for any $n\in\N$ unless $A_n=Y_g$. This choice together with Lemma~\ref{lem:big}(ii) ensures
$$\lim_{n\to\infty}\frac{|A_n\cap \beta(G)|}{|A_n|}\ =0.$$
Therefore, the sequence $\{A_n\}_{n\in \N}$ is a F{\o}lner sequence for the $G$\ti action on $X$ no matter how we define $\beta|_D:D\to H$. It follows (Remark~\ref{rem:gen}) that the $G*H$\ti action on $X$ is amenable. On the other hand, this action is transitive regardless of the definition of $\beta|_D$ since the $G$\ti orbit $\beta(G)$ meets every $H$\ti orbit.

\smallskip

It remains to show that the action is faithful upon suitably determining $\beta$ on $D$. We will refer to elements of $G*H$ written in their reduced canonical form as {\it words}. With a customary abuse of notation, a (non-trivial) word is of the form $w = h_n g_n \ldots h_2 g_2 h_1 g_1$ where only $g_1$ and $h_n$ are allowed to be trivial. Since there are only countably many words, it is enough to prove the following claim: For any $w$ and any finite set $E\se D$ on which $\beta$ is already determined, we can prescribe $\beta$ on a finite set $E'\se D$ such that (i)~$E\cap E'=\varnothing$ and $\beta|_{E\sqcup E'}$ is injective, (ii)~for any injection $\beta|_D:D\to H$ extending $\beta|_{E'}$, there is $x_0\in \beta(E')$ with $wx_0\neq x_0$ (hence $w$ is not in the kernel of $G*H\to \langle G_\beta, H\rangle$).

\smallskip

In order to prove the claim, we assume for definiteness that $g_1, h_n\neq e$ (the other cases are similar). We use indices $1\leq i\leq n$ and $0\leq j\leq n$. Since $H$ is infinite, we may pick $x_0, y_i\in H$ such that all the $2n+1$ elements $x_0,\, y_i,\, h_i y_i,$ are distinct and not in $\beta(E)$. Further we pick $d_j\in D$ such that all the $2n+1$ elements $d_j,\, g_i d_{i-1}$ are distinct and not in $E$. We set $E'=\{d_j, g_i d_{i-1}\}$ and define $x_i=h_i y_i$. The prescriptions $\beta(d_j) = x_j$, $\beta(g_i d_{i-1}) = y_i$ extend $\beta$ injectively; moreover $w x_0 = x_n$ is indeed different from $x_0$. This proves the claim and thus concludes the proof.
\end{proof}

\subsection{}
Another ingredient for Theorem~\ref{thm:main} is:

\begin{theorem}\label{thm:upgrade}
Let $G$, $H$ be countable groups and assume that $G*H$ has a transitive amenable action admitting a F{\o}lner sequence $\{A_n\}_{n\in\N}$ with $|A_n|$ unbounded. Then $G*H\in\A$.
\end{theorem}

\begin{corollary}\label{cor:fin:quot}
Let $G$ and $H$ be finitely generated countable groups both admitting a finite index proper subgroup. Then $G*H\in\A$.
\end{corollary}

\begin{proof}[Proof of the corollary]
Let $G_0\lneqq G$ and $H_0\lneqq H$ be finite index subgroups; we may assume that they are normal. The free product of the quotients $L=(G/G_0)*(H/H_0)$ admits a (non-trivial, finitely generated) free subgroup $F<L$. If $F$ has rank one it is amenable and hence in~$\A$. Otherwise $F$ is of the form $\Z*\cdots*\Z$ and hence is also in~$\A$ by Theorem~\ref{thm:weaker} (or already by van Douwen's result). It follows from Lemmata~\ref{lemma:ind0} and~\ref{lemma:ind} that $L$ is in~$\A$. We thus have a corresponding transitive amenable action $G*H\to L \to X!$ and need only show it satisfies the assumption of Theorem~\ref{thm:upgrade}. If not, there would be a F{\o}lner sequence $\{A_n\}_{n\in\N}$ of cardinality bounded by some $N\in\N$. Let $S$ be a finite set generating $G*H$. In particular, since some $A_n$ is a $(S,\epsi)$\ti F{\o}lner set for $\epsi<1/N$ we obtain a non-empty finite $G*H$\ti invariant subset $A_n\se X$, which is impossible since $X$ is countable and the action transitive.
\end{proof}

\begin{proof}[Proof of Theorem~\ref{thm:upgrade}]
Let $Y$ be a $G*H$\ti set as in the assumption and choose any $y_0\in Y$. Fix two non-trivial elements $g_0\in G$, $h_0\in H$ and consider the transitive $G*H$\ti set $Z=G*H/\langle g_0 h_0\rangle$; let $z_0\in Z$ be the trivial coset $\langle g_0 h_0\rangle$. It follows e.g. from elementary Bass-Serre theory~\cite{Serre77} that (i)~$z_0\neq h_0 z_0 = g_0^{-1} z_0$ and (ii)~any non-trivial $w\in G*H$ fixes at most one point in $Z$. The latter property can be verified as follows. If $w$ has fixed points, we may assume upon conjugating that $w=(g_0 h_0)^n$ for some $n\in\Z$. Now any $w$\ti fixed point $u\langle g_0 h_0\rangle$ satisfies $(g_0 h_0)^n u\langle g_0 h_0\rangle = u\langle g_0 h_0\rangle$. By uniqueness of the normal form for free products, $u$ is in $\langle g_0 h_0\rangle$ and hence determines the same fixed point.

We set now $X=Y\sqcup Z$ and endow it with the corresponding action $G*H\to X!$. Consider the permutation $\sigma\in X!$ that transposes $(y_0,z_0)$ and is trivial otherwise. We claim that the new $G*H$\ti action obtained by conjugating $H$ by $\sigma$ is faithful, transitive and amenable. To be more precise, we consider the original $G$\ti action on $X$ and the new $H$\ti action given by $(h,x)\mapsto h^\sigma x = \sigma^{-1} h \sigma x$ (note $\sigma^{-1}=\sigma$). This yields a new $G*H$\ti action by the universal property of free products; by abuse of notation we denote this action by $(w,x)\mapsto w^\sigma x$ for $w\in G*H$.

\smallskip

\emph{Amenability.} Upon passing to a subsequence, we may assume $|A_n|\to \infty$. But for any $s\in G*H$ there are only finitely many points for which the new action is different from the original one, so that we still have $\lim_{n\to \infty} |A_n \vartriangle s A_n|/|A_n|=0$.

\smallskip

\emph{Transitivity.} We claim that every point lies in the orbit of $y_0$ and will use~(i) above. First pick any $y\in Y$, $y\neq y_0$. Let $w\in G*H$ be the \emph{shortest} word such that $w y_0=y$ in the original action; $w$ exists since $G*H\to Y!$ was transitive. If the rightmost letter of $w$ is in $G$, then $w^\sigma y_0 = w y = y$ (because the ``trail'' of $y_0$ under all non-empty right prefixes of $w^\sigma$ avoids $y_0$ by minimality of $w$). Otherwise, the rightmost letter of $w$ is some $h\in H$ with $h y_0 \neq y_0$. But then $h y_0 = h^\sigma g_0 h_0^\sigma y_0$ and thus replacing $h$ with $h g_0 h_0$ yields a new word $v$ with $v^\sigma y_0 = y$. Now pick any $z\in Z$. If $z=z_0$ we attain it since $z_0 = g_0 h_0^\sigma y_0$. Otherwise, let $w$ be the shortest word such that $w z_0 = z$ in the original action. If the rightmost letter of $w$ is in $G$, then again $w^\sigma z_0 = w z_0 = z$ and hence $z = (wg_0 h_0)^\sigma y_0$. Otherwise this letter is some $h\in H$, in which case $w^\sigma y_0 = z$ since $h^\sigma y_0 = h z_0 \neq z_0$.

\smallskip

\emph{Faithfulness.} Let $w\in G*H$ be any non-trivial word. In view of~(ii) above, there is some $z\in Z$ such that in the original action not only $w z \neq z$ but also the ``trail'' of $z$ under all right prefixes of $w$ avoids $z_0$. Therefore $w^\sigma z = w z \neq z$ and the action is faithful.
\end{proof}

\subsection{}
We are now ready to finish the

\begin{proof}[Proof of Theorem~\ref{thm:main}]
In view of Theorem~\ref{thm:weaker}, we may assume without loss of generality that both $G$ and $H$ have finite index subgroups $G_0\lneqq G$ and $H_0\lneqq H$ with the fixed point property~(F). By Lemma~\ref{lem:K}, a non-finitely-generated countable group cannot have virtually~(F). Therefore, both $G$ and $H$ are finitely generated and we conclude from Corollary~\ref{cor:fin:quot} that $G*H\in\A$ in this case too. The converse was established as Lemma~\ref{lem:converse}.
\end{proof}

\begin{proof}[Proof of Corollary~\ref{cor:infinite}]
Keep the notation of the corollary. If $n=\infty$, then $H=*_{i=2}^\infty G_i$ is not finitely generated and in particular by Lemma~\ref{lem:K} does not virtually have~(F). Thus we may apply Theorem~\ref{thm:main} and conclude that $G=G_1*H$ is in~$\A$.

Assume now $n\neq\infty$. If all $G_i$ with $1<i\leq n$ have property~(F) and $G_1$ virtually has~(F), then Lemma~\ref{lem:converse} and Lemma~\ref{lem:F:products} show that $G=G_1*(G_2 * \cdots *G_n)$ has virtually~(F). If not, we may assume upon reordering that either (i)~both $G_1$ and $G_2$ fail to have property~(F) or (ii)~$G_1$ does not virtually have~(F). In case~(i), $G_2* \cdots * G_n$ also fails to have property~(F) and hence $G=G_1*(G_2 * \cdots *G_n)$ is in~$\A$ by Theorem~\ref{thm:main}. In case~(ii) Theorem~\ref{thm:main} implies $G\in\A$ aswell.
\end{proof}

\section{Remaining proofs}
\label{sec:little}
\subsection{}
In this section, all Roman numerals refer to the properties listed in Theorem~\ref{thm:properties}. Point~(\ref{P:amenable}) follows by definition. Point~(\ref{P:*}) follows from Theorem~\ref{thm:main}. For~(\ref{P:product}), let $G\to X!$ and $H\to Y!$ be faithful transitive amenable actions. Then the $G\times H$\ti action on $X\times Y$ is faithful and transitive. If $A_n\se X$ and $B_n\se Y$ form F{\o}lner sequences, then $A_n\times B_n$ yields a F{\o}lner sequence for the product. Conversely, if $G\times H\to X!$ is amenable, then, as a $G$\ti set, $X$ is isomorphic to $X_0\times X_1$, where $X_0$ is any $G$\ti orbit in $X$ and $X_1=G\bsl X$ with trivial $G$\ti action. The $G$\ti action is faithful, transitive and amenable by Lemma~\ref{lemma:map} since there is a (non-canonical) $G$\ti map $X\to X_0$.

\smallskip

The first part of~(\ref{P:co-amenable}) follows from Lemmata~\ref{lemma:ind0} and~\ref{lemma:ind}. For the second part, let $Q$ be any group and consider $G=H\ltimes {\bf Z}$ where $H=\bigoplus_{\bf Z} Q$ denotes the direct sum (=restricted product) and ${\bf Z}$ acts on it by shift. It was observed in~\cite{Monod-Popa} that the subgroup $K=\bigoplus_{\bf N}Q$ is co-amenable in $G$. On the other hand, the $G$\ti action on $X=G/K$ is faithful and transitive; thus $G\in\A$. Since $H\cong Q\times H$, it is enough now by~(\ref{P:product}) to chose for $Q$ any group not in~$\A$ (e.g. using Lemma~\ref{lem:T} below).

\subsection{}
For the sake of our discussion, it is convenient to introduce the following.

\begin{definition}
Let~$\B$ be the class of all countable groups admitting some amenable action on a countable set without finite orbits.
\end{definition}

Whilst $\B$ contains $\A$, it is much wider; notice for instance that any group with a quotient in~$\B$ is itself in~$\B$. Moreover (Lemma~\ref{lem:K}), Kazhdan's observation~\cite{Kazhdan67} shows that any countable group that is not finitely generated is in~$\B$. For example, let $G=\bigoplus_{n=1}^\infty G_n$, where $G_n$ are any countable groups with $G_1\notin\A$; then $G\notin \A$ by~(\ref{P:product}) but $G\in\B$.

A group in~$\B$ cannot have virtually~(F); however, the class of groups to which Theorem~\ref{thm:main} applies is yet much wider than~$\B$, since it contains notably every countable group without minimal finite index subgroup (e.g. $\mathbf{SL}_3(\Z)$). In summary:
$$\A\ \subsetneqq\ \B\ \subsetneqq\ \big\{ \text{ not virtually~(F) }\big\}.$$

\subsection{}
Another well known criterion for amenability is that the $G$\ti action on $X$ is amenable if and only if the associated unitary representation on $\ell^2(X)$ \emph{almost has invariant vectors} in the sense of Kazhdan~\cite{Kazhdan67}; see also~\cite{Harpe-Valette} for this notion, Kazhdan's property~(T) and the Kazhdan-Margulis relative property~(T). 

Since $\ell^2$\ti functions have finite level-sets, it follows:

\begin{lemma}\label{lem:T}
Every amenable action of a group with Kazhdan's property~(T) has a finite orbit. In particular, any Kazhdan group without non-trivial finite quotients has property~(F).\hfill\qedsymbol
\end{lemma}

\noindent
Thus a Kazhdan group can never be in~$\B$ and~(\ref{P:T}) follows.

\smallskip

Notice that the second statement of Lemma~\ref{lem:T} applies in particular to infinite simple Kazhdan groups; we point out that such groups do indeed exist, as follows from the work of Gromov~\cite{Gromov87} and Ol{\cprime}shanski\u{\i}~\cite{Olshanskii93} on quotients of hyperbolic groups (see also~\cite{Ozawa04}).

\smallskip
There exist however groups with property~(F) that are not Kazhdan groups; indeed, property~(F) is stable under free product (Lemma~\ref{lem:F:products}) whilst a (non-trivial) free product is never a Kazhdan group~\cite{Harpe-Valette}. Since property~(F) passes to quotients, we remark further that any group generated by finitely many subgroups with property~(F) still enjoys this property. (See also Lemma~\ref{lem:monster} below.)

\subsection{}
Point~(\ref{P:rel_T}) follows from an argument similar to~(\ref{P:T}):

\begin{lemma}
Let $H\triangleleft G$ be a normal subgroup of $G\in\A$. If the pair $(H,G)$ has the relative property~(T), then $H$ has finite exponent.
\end{lemma}

\begin{proof}
Let $G\to X!$ be a faithful transitive amenable action. Since $G$ almost has invariant vectors in $\ell^2(X)$, the relative property~(T) implies that $H$ fixes a non-zero vector in $\ell^2(X)$. Therefore $H$ preserves a non-empty finite set in $F\se X$. Since $H$ is normal in $G$ and $G$ acts transitively, it follows that $X$ has a partition into $G$\ti translates $gF$ of $F$ and that this partition is preserved by $H$. Therefore we have a natural morphism $H\to \prod_{gH\in G/H} (gF)!$. Since the $G$\ti action is faithful, this morphism is injective and thus realises $H$ as a subgroup of a group of exponent $|F|!$.
\end{proof}

Now~(\ref{P:semidirect}) also follows: Since $\mathbf{SL}_2({\bf Z})$ is virtually free, it is in~$\A$ by~(\ref{P:co-amenable}); the group ${\bf Z}^2$ is in~$\A$ since it is amenable. On the other hand, the natural semi-direct product ${\bf Z}^2\rtimes\mathbf{SL}_2({\bf Z})$ is not in~$\A$ since it has the relative property~(T) (see~\cite{Harpe-Valette}). This is yet another example of a group in $\B\setminus\A$, since its quotient $\mathbf{SL}_2({\bf Z})$ is in~$\A$.

\subsection{}
At this point we have a number of alternative proofs of~(\ref{P:embed}): Any countable group $Q$ embeds into some $G\in\A$. For instance, take $G=Q*{\bf Z}$, or $G=*_{n=1}^\infty Q$ (by Corollary~\ref{cor:infinite}), or the group $G=H\rtimes{\bf Z}$ with $H=\bigoplus_{\bf Z}Q$ considered earlier. (See also~\cite{Grigorchuk-Nekrashevych}.)

\subsection{}
It is unknown whether R.~Thompson's group $F$ is amenable or not. This group is defined in detail in~\cite{Cannon-Floyd-Parry}; all we need to know here is that it satisfies the assumptions of the following:

\begin{lemma}\label{lem:F}
Let $F$ be a group of orientation-preserving piecewise linear homeomorphisms of the interval $(0,1)$. If the derived subgroup $F'$ has a dense orbit, then $F'$ and $F$ are in~$\A$.
\end{lemma}

\begin{proof}
Let $X=F'x$ be such an orbit and notice that $F'$ acts faithfully on $X$. Choose a sequence $\{x_n\}$ in $X$ converging to~$0$. Let $\mu$ be a limit point of the sequence of point-measures $\delta_{x_n}$ in the space of means on $X$. Then $\mu$ is $F'$\ti invariant because for every $g\in F'$ there is $\epsi>0$ such that $g$ is trivial on $(0, \epsi)$ (this is where we use the assumptions on $F$). Therefore $F'\in\A$ and thus $F\in\A$ by~(\ref{P:co-amenable}).
\end{proof}

\subsection{}
A \emph{Tarski monster} is a non-cyclic group such that all its proper subgroups are cyclic. Ol{\cprime}shanski\u{\i} constructed various Tarski monsters and then proved that his groups are non-amenable~\cite{Olshanskii80}. Therefore the following applies to them:

\begin{lemma}\label{lem:monster}
Let $G$ be a non-amenable group such that all its proper subgroups are amenable. Then $G$ has property~(F).
\end{lemma}

\begin{proof}
Suppose for a contradiction that $G$ has an amenable action on the countable set $X$ without fixed points. Then the stabiliser of any point is amenable and it follows that the action is amenable in Zimmer's sense~\cite{Zimmer84}. However, a group $G$ is amenable if and only if some (or equivalently all) of its actions on a countable set are both amenable and amenable in Zimmer's sense. For the reader's convenience, we sketch the argument without appealing to Zimmer's notion: For any finite set $S\se G$ and $\epsi>0$ we may choose an orbit $G/L\se X$ and an $(S,\epsi/2)$\ti F{\o}lner set $A\se G/L$. Fix a section $\sigma:G/L\to G$ of the natural map $G\to G/L$. Since $L$ is amenable, there is a $(T,\epsi/2)$\ti F{\o}lner set $B\se L$ for $T=\{\sigma(sgL)^{-1}s\sigma(gL) : s\in S, gL\in A\}\se L$. One verifies that $\{\sigma(gL)\ell: gL\in A, \ell\in B\}$ is an $(S,\epsi)$\ti F{\o}lner set in $G$. Thus $G$ is amenable.
\end{proof}

\subsection{}

Finally we propose a problem: What can be said about amenable actions of the Hilbert modular group $\Gamma=\mathbf{SL}_2({\bf Z[\sqrt{2}]})$~? In particular, does $\Gamma$ belong to any of the classes considered in this note?

More generally, let $\Gamma$ be a lattice in a product $G=G_1\times G_2$ of locally compact groups and assume that both projections $\Gamma\to G_i$ have dense image. Are there natural conditions on the groups $G_i$ that imply $\Gamma\in\A$ or $\Gamma\notin\A$~? (Not desired here are conditions so coarse as to apply to \emph{any} lattice in $G$, such as a product of lattices $\Gamma_i<G_i$.)

\smallskip

In another direction, are there interesting amenable actions of lattices in $\mathbf{SO}(n,1)$ or in $\mathbf{SU}(n,1)$~?


\end{document}